\documentclass[12pt]{extarticle}
\usepackage[margin=1in]{geometry}
\usepackage{dsfont}
\usepackage{amsfonts}
\newtheorem{theorem}{Theorem}[section]
\newtheorem{pro}{Proposition}[section]
\newtheorem{lemma}{Lemma}[section]
\newtheorem{definition}{Definition}[section]
\newtheorem{remark}{Remark}[section]
\newtheorem{cor}{Corollary}[section]

\newcommand{\proof}[1]{\noindent{\it\bf Proof:#1\ }}
\newcommand{\QED}{\hfill$\Box$\medskip}

\begin{document}

\title{   On   Action of  $PSL(n+1, {\bf C})$ on Space of  $L_k^p$-maps on ${\bf P}^n$ }
\author{Gang Liu\\Department of Math, UCLA  }
\date{ October,  2018}
\maketitle

\begin{abstract}
 In  this paper, we  prove properness  of the action of the reparametrization group $PSL(n+1, {\bf C})$ on  the space of $v$-stable $L_k^p$-maps on ${\bf P}^n$  as well as  related results. They   extend    our  earlier work on the proper action of  
the reparametrization groups on the space of weakly stable nodal $L_k^p$-maps.

\end{abstract}
 
\begin{large}
\section{Introduction}

In  \cite{Liu2, Liu3} we have introduced the notion of weak-stable nodal $L_k^p$-maps
 as a natural generalization of  stable $J$-holomorphic maps in Gromov-Witten theory introduced by Kontsevich in \cite{Kontsevich}.
The purpose of this paper  is to generalize part of the results
 on proper  action of the reparametrization  groups on  the space of stable or weakly stable $L_k^p$-maps \cite{Liu2, Liu3} to the case  with  domain  ${\bf P}^n$ acting on by $PSL(n+1, {\bf C})$. Further generalizations  to  other 
higher dimensional but smooth domains will be treated in a companion of this paper.

The main difficulty for such a generalization lies on the well-known fact:  unlike $1$-dimensional case, higher dimensional  biholomorphic maps are not conformal in general. Thus even though a version of   higher dimensional  energy function  of  the  $L_k^p$-maps,  the key quantity  used  for the $1$-dimensional  case \cite{ Liu3},
can be  defined, it is only invariant with respect to the  conformal  automorphisms of the domains but not respect to  the natural action of the  reparametrization  group of  biholomorphic automorphisms.

In this paper, instead of using energy function,  the  volume function $v(f)$   is used to define 
the notion of $v$-stability for $L_k^p$-maps and to prove part of the corresponding results in \cite {Liu3} described as follows.

Let $M$ and  $N$ be   compact  Riemannian manifolds of class $C^{\infty}$ with $\dim (M)\leq \dim (N)$  and ${{\cal M}}_{k,p}(M, N)$ be the mapping space  of  $L_k^p$-maps.
Throughout the paper, we always assume that $m_0=[k-\frac{m}{p}]\geq 1$ where $m=dim
(M)$ so that each $L_k^p$-map is at least of class $C^1$ by Sobolev embedding theorem.

When the domain $M$ is ${\bf P}^n$, the group $G=PSL(n+1, {\bf C})$  operates on  the space ${{\cal M}}_{k,p}({\bf P}^n, N)$ as the group of  reparametrizations.

Recall in the case that $X$ is a locally compact topological space such as a finite dimensional manifold and $G$ is a Lie group or a locally compact topological group, a (continuous) group action $\Phi:G\times X\rightarrow X$ is said to be proper if the map  $\Phi\times{\it id}_{X}:G\times X\rightarrow X\times X$  is a proper map:  the inverse image $(\Phi\times{\it id}_{X})^{-1}(K) $ of any compact subset $K\subset X\times X$ is compact.

In our infinite dimensional case, ${{\cal M}}_{k,p}({\bf P}^n, N)$
is not locally compact, the above definition is too weak to be useful.

In \cite {Liu2, Liu3}, we have introduced the following  stronger definition.
\begin{definition}
	A group action $\Phi:G\times X\rightarrow X$ is said to be proper if  for any compact subset $K\subset X\times X$, there is a neighborhood $U$ of $K$ in  $ X\times X$ such that the image  $\pi_G ((\Phi\times{\it id}_{X})^{-1}(U))$ of the projection to $G$ of the inverse image $(\Phi\times{\it id}_{X})^{-1}(U)$  is pre-compact in $G$.

\end{definition}
It was proved in \cite {Liu2, Liu3} that in the case $X$ is a locally compact topological space acting by a  locally compact  group $G$, the definition here is equivalent to the usual one.

Now we  define  the notion of $v$-stability \cite{Liu2}.

\begin{definition}
	
	A $L_k^p$-map $f:M\rightarrow N$ with  $dim (M)=m$ is said be to $v$-stable if its volume $v(f)>0$. Here $v(f)= v_m(f)$ if $v_m(f)>0$ and $v(f)= v_{m-1}(f)$ if $v_m(f)=0$,  defined by using   the  two  volume functions $v_{m-1}$ and $v_m$ of dimension $m-1$ and $m$ respectively.
	
\end{definition}

The notion of $v_{m-1}$/ $v_{m}$-stability is defined similarly.

\begin{remark}
In the case that $dim (M)=2$, $f:M\rightarrow N$ is $v$-stable  if  and only if $f$ is not a constant map,  and hence weakly stable by the definition in   \cite {Liu3}.  In this sense $v$-stability  is the natural generalization of the notion of weak stability in  \cite {Liu3} to the higher dimensional case with smooth domains. The theorems below show that $v$-stable maps have  the similar properties as weakly stable maps.  
\end{remark}
 
 Let ${{\cal M}}^{*}_{k,p}({\bf P}^n, N)$ be the space of $v$-stable
$L_k^p$-maps on ${\bf P}^n$.

\begin{theorem}
	
The action of $G=PSL(n+1, {\bf C})$ on 	the space ${{\cal M}}^{*}_{k,p}({\bf P}^n, N)$  is  proper.
\end{theorem}

 \begin{cor} 
 	For any    $v$-stable $L_k^p$-map $f$,  its stabilizer $\Gamma_f$ is always a compact
 	subgroup of $G$.
 \end{cor}

\begin{pro}
	
	The action of $G=PSL(n+1, {\bf C})$ on 	the space ${{\cal M}}^*_{k,p}({\bf P}^n, N)$    is $G$-Hausdorff. Therefore, the  quotient space  ${{\cal M}}^*_{k,p}({\bf P}^n, N)/G$  of unparametrized $v$-stable  $L_k^p$-maps is Hausdroff.
\end{pro}

\begin{cor} 
	Given any    
	$f$ in ${ {\cal M}}^*_{k, p}$, the G-orbit $G\cdot f$ is closed in ${\cal M}^*_{k, p}$.  
\end{cor}

\begin{definition} 
	A $L_k^p$-map $f$ in ${ {\cal M}}^*_{k, p}$ is  said to be  stable if	  its stabilizer $\Gamma_f$ is  a  finite group.
	
\end{definition}

\begin{theorem}
	Any  $v_{2n}$-stable $L_k^p$-map on ${\bf P}^n$ is stable.	
\end{theorem}

 \begin{definition}
 A $L_k^p$-map $f:{\bf P}^n\rightarrow N$ is said to be (one of ) the standard $S^1$-invariant map if $f={\bar f} 	\circ \pi$ where $\pi:{\bf P}^n\rightarrow {\bf P}^n/S^1$ is the quotient map of the standard $S^1$-actions on ${\bf P}^n$ and ${\bar f}:{\bf P}^n/S^1\rightarrow N$ is the induced map.
 
 \end{definition}

\begin{theorem}
	Any  $v$-stable $L_k^p$-map on ${\bf P}^n$ is stable if it is not (one of ) the standard $S^1$-invariant map (up to a conjugation).	
\end{theorem}

\begin{cor} 
	 For a $v$-stable $L_k^p$-map, if the identity component $\Gamma_f^0$ of the isotropy group   is nontrivial, it is isomorphic to $S^1$. This can happen only when $f$ is one of  the standard $S^1$-invariant map (up to a conjugation).
	
\end{cor}

After  defining  the two volume functions in Sec. 2, the theorems stated  in this section are proved in Sec. 3.

\section{ Definitions of the Volume Functions  }

In this section we give the definitions $v_m(f)$ and $v_{m-1}(f)$ of  the volume functions.  
 
 The definition for $v_m(f)$ is standard that we recall now.
For  Reimannian manifolds $M$ and $N$ and a $C^1$ map $f:M\rightarrow N$ with $dim(M)=m$,  the $m$-dimensional volume function $v_m(f):=vol_m(f)=\int_M v_m(df) d\nu_M,$ where $d\nu_M=|d{\tilde \nu_M}|$  is the volume density  and  $d{\tilde \nu_M} $  is the  volume  form  determined by the  metric of $M$. Here $d{\tilde \nu_M} $ is defined upto a sign and 
$v_m(df)=:vol_m(df)$ is a non-negative function defined on $M$
as follows. For any $x\in M,$ let $(e_1, \cdots, e_m)$ be a orthonormal frame  of $T_xM$. Then $v_m(df)(x)$ is defined to be the  $m$-dimensional volume of the parallelepiped spanned by $df (e_1), \cdots, df (e_m)$ in $T_{f(x)}N$ measured by the metric $g_{f(x)}$ on $T_{f(x)}N.$ Note that $v(df)(x)$ is independent of the choices of the orthonormal frames  of  $T_xM$ since $|df (e_1)\wedge  \cdots\wedge df (e_m)|=|det (a_{i, j}) df (e'_1)\wedge  \cdots\wedge df (e'_m)|$ where $ e_i=\Sigma_{j}a_{i, j} e'_j$.

In the case that $f$ is a immersion at $x$, hence a local embedding on a small neighborhood $U$ of $x$ with the $m$-dimensional image ${\tilde U}$,  consider the restriction map $f:U\rightarrow {\tilde U}\subset N$ and let $d{\tilde \nu}_{m,\tilde U}$ be the $m$-dimensional volume form  (defined upto a sign) with respect to the induced metric.
On $U$, define 	${\tilde v}_m(df)=f^*(d{\tilde \nu}_{m, \tilde U})$. Then   $v_m(df)=|{\tilde v}_m(df)/d{\tilde \nu}_{m, U}| $.

This implies the next  lemma.

\begin{lemma} 
	Let $u:M_1\rightarrow M$ be a diffeomorphism 
	and $ f:M\rightarrow N$ be a $C^1$-map. 
	Then  in a small neighborhood of  any $x\in M_1$ where  $f\circ u$ is a immersion, ${\tilde v}_m(d(  f\circ u))= u^*[{\tilde v}_m(d( f))].$

\end{lemma}

\begin{cor}
	Let $R_1$ be a finite region in $M_1$ in the above lemma  then $v_m(f\circ u|_{R_1})=v_m(f|_{u(R_1)})$.
	
\end{cor}

\proof

Let $M^*$ be the open subset of $M$ such that at any point $x\in M^*$ $f:M\rightarrow N$ is a immersion  and $M_1^*= u^{-1} (M^*)$ be the corresponding open set in $M_1$. Clearly $v_m(f|_{M^*})=v_m(f)$ and $v_m(f\circ u^{-1}|_{M^*_1})=v_m(f\circ u^{-1})$ by the definition.
Hence we may assume that $R_1$ is lying inside $M^*_1$ and $u(R_1)$ is lying inside $M^*.$

By using a partition of unit, we only need to consider the local case where we may assume that everything involved is oriented and orientation preserving.  
 Then  $v_m(df)d\nu_{M}={\tilde v}_m(df).$ 

By the lemma, $$v(f\circ u|_{R_1})=\int_{R_1} {\tilde v}_m(d(f\circ u))=\int_{R_1} u^*({\tilde v}_m(df))
$$ $$ =\int_{u(R_1)}(u^{-1})^* [u^*({\tilde v}_m(df))]
=\int_{u(R_1)}{\tilde v}_m(df)=v(f|_{u(R_1)}).$$

\QED

\begin{definition}
	A $L_k^p$-map $f:M\rightarrow N$ is said to be $v_m$-stable if $v_m(f)>0$.
\end{definition}

 The $(m-1)$-volume function $v_{m-1}(f)$ is only defined for
 $f:M\rightarrow N$ that is not  $v_m$-stable. For such an $f$, 
$v_m(df)(x)=0$ or equivalently $rk (df_x)<m$  for any $x\in M$. 

In this situation, let $M_{m-1}^*(f)$ be the open subset of $M$ consisting of points $x$ where $rk (df_x)=m-1$. If $M_{m-1}^*(f)$ is empty, define $v_{m-1}(f)=0.$ Otherwise $v_{m-1}(f)$ is defined as follows.

For any point $x\in M_{m-1}^*(f)$, by definition the kernel $K_x=\ker (df_x)$ is $1$-dimensional subspace of the tangent space $T_xM$. The collection $K$ of $K_x$ with $x\in M_{m-1}^*(f)$ is an  $1$-dimensional
 distribution on  $M_{m-1}^*(f)$ and  hence is  integrable. Thus we get an $1$-dimensional foliation on $M_{m-1}^*(f)$. Locally  near a given point
 $x_0\in M_{m-1}^*(f)$ by choosing a local  $C^0$-section  of $K$ of unit length, we get an ordinary differential equation so  the foliation near $x_0$ can be ``represented'' by the integral curves of the differential equation. The  differential equation here is  canonically
defined up to a choice of a local orientation for $K$, hence is globally defined on the double covering of $M_{m-1}^*(f)$ even $K$ is not orientable. Let $S(x_0)$ be a local slice at $x_0$ transversal  to the foliation with $S(x_0)\simeq B^{m-1}$ of a small $(m-1)$-dimensional  open ball, and $U_\epsilon(S(x_0) ) \simeq S(x_0)\times (-\epsilon,\epsilon)$ be the collection of integral curves with initial values at $t=0$ lying on $S(x_0)$ and $t$ varying in the sufficiently small interval $(-\epsilon,\epsilon)$. Each such neighborhood $U_\epsilon(S(x_0) )$ will be called a $K$-neighborhood of $x_0$. Then  $M_{m-1}^*(f)$ can be covered by countably many such open $K$-neighborhoods  $U_{\epsilon_i}(S(x_i) )$.

Now fix such a $U_\epsilon(S(x_0) )$ and consider the restriction map $f:U_\epsilon(S(x_0) )\rightarrow N$. Since $f$ is constant along any leaf of the foliation,   the image of $U_\epsilon(S(x_0) )$ is the same as
that for $f:S(x_0) \rightarrow N$, which is a local embedding when $S(x_0)$ is sufficiently small. Let ${\tilde S}(x_0)$ be the image of $S(x_0)\simeq B^{m-1}$ as a $(m-1)$-dimensional submanifold in $N$.
Note that for any other local transversal slice $S(x'_0)$ in  $U_\epsilon(S(x_0) )$, the integral curves induce an identification
$\gamma_{x_0', x_0}:S(x_0)\rightarrow S(x'_0)$ of class at least $C^1$  such that $f|_{S(x'_0)}\circ \gamma_{x_0', x_0}=f|_{S(x_0)}$.
Thus the images ${\tilde S}(x_0)$ and ${\tilde S}(x'_0)$ are the same.

For any point $x\in S(x_0),$ let  $d{\tilde \nu}_{m-1,{\tilde S}(x_0)}$ be the $(m-1)$-dimensional volume form  (defined up to a sign) with respect to the induced metric.

\begin{definition}
On $ S(x_0)$, define 	${\tilde v}_{m-1}(df_{ S(x_0)})=f|_{ S(x_0)}^*(d{\tilde \nu}_{m-1, {\tilde S}(x_0)})$ 

and  $v_{m-1}(df|_{S(x_0)})=|{\tilde v}_{m-1}(df|_{ S(x_0)})/d{\tilde \nu}_{m-1, { S}(x_0)}| $.
\end{definition}

Then  the next  lemma follows from the definitions above.

\begin{lemma} 

	Let $u:M_1\rightarrow M$ be a diffeomorphism 
	and $ f:M\rightarrow N$ as above. 
	Consider a small  transversal slice  $S(y_0)$  in    $( M^*_1)_{m-1}(f\circ u)$. Let   $S(x_0)=u( S(y_0))$ be the corresponding  transversal slice in $M^*_{m-1}(f)$. Then ${\tilde v}_{m-1}(d(  f\circ u|_{S(y_0)}))= u^*[{\tilde v}_{m-1}(d( f|_{S(x_0)}))].$
	
	For any two small local slice  $ S(x_0)$ and  $ S(x'_0)$ in $U_{\epsilon}(  S(x_0))$, ${\tilde v}_{m-1}(d(  f|_{S(x_0)}))= \gamma_{x_0', x_0}^*[{\tilde v}_{m-1}(d( f|_{S(x'_0)}))].$
\end{lemma}

\begin{definition}
	The local $(m-1)$-volume  of $f$ over $U_{\epsilon}(  S(x_0))$ is defined to be $v_{m-1}(f|_{U_{\epsilon}(  S(x_0))})=\int_{S(x_0)} { v}_{m-1}(d(  f|_{S(x_0)}))d{\nu}_{S(x_0)}$, where $d{\nu}_{S(x_0)}$ is the volume  density on $S(x_0)$.
\end{definition}
  
 Then  by above lemma, we have
 
 \begin{cor}
 The local $(m-1)$-volume  $v_{m-1}(f|_{U_{\epsilon}(  S(x_0))})$ is well-defined independent of the choice of the local transversal slice.	
 	Let $u:M_1\rightarrow M$ be a diffeomorphism 
 and $ f:M\rightarrow N$ as above lemma.  Consider a  small $K$-neighborhood  $U(S(y_0))$  in    $ (M^*_1)_{m-1}(f\circ u)$. Let   $U(S(x_0))=u(U( S(y_0)))$ be the corresponding   $K$-neighborhood  $U(S(y_0))$  in   $M^*_{m-1}(f)$. Then ${ v}_{m-1}(  f\circ u|_{U(S(y_0))})=  v_{m-1}(f|_{U(S(x_0))}).$	
 \end{cor}
 
 The proof of the corollary is essentially the same as the proof of the Corollary 2.1. We leave it to the readers.
 
\vspace{2mm}
\noindent 

To define the global $(m-1)$-volume $v_{m-1}(f)$, we define 

\noindent $v_{m-1}(f|_{\cup_{i=1}^lU(S(x_i))})$ for the finite union of $K$-sets $\cup_{i=1}^l U(S(x_i))$ first.  The image $\cup_{i=1}^l{\tilde S}(x_i))$ of $\cup_{i=1}^lU(S(x_i))$ under $f$ is an immersed $(m-1)$ dimensional submanifold with self-intersections.
To define $v_{m-1}(f|_{\cup_{i=1}^lU(S(x_i))})$, we  construct  partitions of unit on subsets of the image $\cup_{i=1}^l{\tilde S}(x_i))$ first. 

Let  ${ {\beta}}_i$ be  a smooth bump-off function supported in $U(S(x_i))$ and $\beta_{S, i}$ be the restriction of $\beta_i$ to $S(x_i)$. Then $\beta_{S, i}$ can be consiered as a function define on the image ${\tilde S }(x_i)$. Consider the open subset $S'(x_i)$ of $S(x_i)$ consisting of points $x$ where ${ {\beta}}_i(x)>0$. Then on the finite union $\cup_{i=1}^l {\tilde S'}(x_i))$ of the images of $S'(x_i),$ the nonzero bump-off functions 
$\beta_{S, i}, i=1, \cdots , l$ define a partition of unit in the usual manner. Denote the bump-off function of the partition of unit on ${\tilde S'}(x_i)\simeq S'(x_i)$ by $\alpha_i$.   

For each choice
 of $\beta$ above, define $$v_{m-1, \beta}(f|_{\cup_{i=1}^lU(S(x_i))})=\Sigma_{i=1}^l|\int_{S'(x_i)} f^*(\alpha_id{\tilde \nu}_{{{\tilde S'}(x_i)}})|$$ $$=\Sigma_{i=1}^l|\int_{S'(x_i)}f^*(\alpha_i) \cdot {\tilde v}_{m-1}(df|_{S'(x_i)})|. $$

\noindent
\begin{definition}
The $(m-1)$-volume functions are   defined by: 

(1) $v_{m-1}(f|_{\cup_{i=1}^l U(S(x_i))})=\sup_{\beta} v_{m-1, \beta}(f|_{\cup_{i=1}^l U(S(x_i))})$ over all bump-off functions $\beta$; (2) $v_{m-1}(f)$ is equal to  the supremum 
of $v_{m-1}(f|_{\cup_{i=1}^l U(S(x_i))})$  over all finite union  of $K$-subsets $ \cup_{i=1}^l U(S(x_i))$ in $M^*_{m-1}.$
\end{definition}

 For the proofs in next section only the properties for the local volume $v_{m-1}$ stated in above lemma are used.  The definition above for the global $(m-1)$-volume is to ensure that if $v_{m-1}(f)>0$,  so is $v_{m-1}(f|_{U(S(x_0))})>0$ for some $K$-set $U(S(x_0))$ in $M^*_{m-1}.$
 
 The functions  $v_i(f)$ for all  $0<i\leq m=\dim(M)$ can be defined using similar ideas. However for the proofs in this paper only above two functions are useful.
 
\section{  Proof  of the Main Theorems }

 \subsection{  Proof  of Theorem 1.1, Proposition 1.1, Corollary 1.1 and Corollary 1.2}
 \vspace{2mm}
 \noindent 
 
We make a reduction first.

\begin{lemma}
	A group action $\Phi:G\times X\rightarrow X$ is proper	if and only if for
	any point $p\in X\times X$ there is a neighborhood $U$ of $p$ in  $ X\times X$
	such that
	the closure of  $\pi_G ((\Phi\times{\it id}_{X})^{-1}(U))$ is compact in $G$.
	
\end{lemma}

The proof of the lemma is elementary and is given in \cite{Liu3}.

By the above lemma,  the properness of the action of $G$ on  ${{\cal M}}^{*}_{k,p}({\bf P}^n, N)$   can be derived from the following theorem.

\begin{theorem}
	The  action of $G={PSL(n+1, {\bf C})}$ on ${{\cal M}}^{*}_{k,p}({\bf P}^n, N)$  has the following property:
	for any $f_1$ and $f_2$  there exist the   open neighborhoods $U_{\epsilon_1}(f_1)$ and $U_{\epsilon_2}(f_2)$ containing $f_1$ and $f_2$   and compact subsets $K_1$ and $K_2$ in 	$G$  accordingly   such that  for any $h_1$ in $U_{\epsilon_1}(f_1)$ ($h_2$ in $U_{\epsilon_2}(f_2)$) and $g_1$ in $G\setminus K_1$  ($g_2$ in $G\setminus K_2$ ), $g_1\cdot h_1$ is not in $U_{\epsilon_2}(f_2)$ ($g_2\cdot h_2$ is not in $U_{\epsilon_1}(f_1)$).

\end{theorem}

\proof

We start with some elementary linear algebra. 
For any  $g \in SL(n+1, {\bf C}),$ 
we have a decomposition in  $SL(n+1, {\bf C}),$ $g=h\cdot  u$ with $u\in SU(n+1)$ and $h$ being self-adjoint and positive.  Indeed  $h=(g\cdot g^*)^{\frac {1}{2}}\in SL(n+1, {\bf C})$ and $u= (g\cdot g^*)^{-\frac {1}{2}}\cdot g\in SU(n+1)$.
Consider the  decomposition $h=w^*\cdot diag ( r_1, r_2\cdots,  r_{n+1})\cdot w$.  Then $g=w^*\cdot diag ( r_1, r_2\cdots,  r_{n+1})\cdot wu.$  Here $r_i>0$ for $ i=1, \cdots , n+1.$ Rename $w^*$ as $u$ and $wu$ as $v$. Denote $diag ( r_1, r_2\cdots,  r_{n+1})$ by $\Delta({\bf r})$ for short. Then we have the decomposition $g=u\cdot \Delta({\bf r})\cdot v$ in $SL(n+1)$ with $u$ and $v$ in 
$SU(n+1)$. This decomposition is not unique for a non generic $g$, but we only need the existence of the decomposition.

We always assume that ${\bf r}$ is ordered  as $0<r_1\leq r_2\leq \cdots \leq r_{n+1}.$ Then as an element in $PSL(n+1)$,  we may assume that $r_{n+1}=1$ and $\Delta({\bf r})=diag (r_1, r_2\cdots,  r_{n}, 1).$ Denote the smallest element $r_{1}$  by $a$  and $ \Delta({\bf r})$
by    $\Delta(a)$.

Assume that the Theorem 2.1 (a)  is not true. Then for any neighbourhoods
$U_{\epsilon_i}(f_i), i=1, 2 $ and any nested  sequences of compact sets $K_1\subset K_2\subset \cdots \subset K_l\cdots$ in $G$, there  are  sequences
$\{g_k\}_{k=1}^{\infty}$ in $G$ and $\{h_k\}_{k=1}^{\infty}$ in $U_{\epsilon_1}(f_1)$ such that (a) $g_k$ is not in $K_{i(k)}$; (b) $h_k\circ
g_k$ is in  $U_{\epsilon_2}(f_2).$

Here $U_{\epsilon_1}(f_1)$,  $U_{\epsilon_2}(f_2)$ and   $K_k, k=1, \cdots$, 
will be decided below. Note that  we allow $f_1=f_2$ but the choice of  $U_{\epsilon_1}(f_1)$ and   $U_{\epsilon_2}(f_2)$ below are  different.
\noindent  

Let $D_{n+1}$ be the collection of all
non-singular  diagonal matrices  with $n+1$ positive  entries.
Choose ${\tilde K}_k\subset  SU (n+1)\times D_{n+1} \times SU (n+1)$  to be $\{(u, \Delta({\bf r}), v)\in SU (n+1)\times D_{n+1} \times SU (n+1)\,|\,\frac {1}{n} \leq r_i\leq 1, i=1, \cdots, k+1\}$, where $\Delta({\bf r})=diag (r_1, \cdots,
r_{n+1})$.  Denote  the corresponding compact set in  $PSL(n+1)$ by ${ K}_n$ obtained by sending $(u, {\bf r}, v)$ to $u\cdot \Delta({\bf r})\cdot v$. 

 First fix  $U_{\epsilon_1}(f_1)$ without any restrictions.
 
 Now we  need to deal with the two case:  (I) $v_m(f_2)>0$ and   (II) $v_m(f_2)=0 $ but $v_{m-1}(f_2)>0$. Here and below $m=2n$.

 For case (I) we may assume that  $v_m(f_1)>0$ as well. Indeed if in this case $v_m(f_1)=0$ but $v_{m-1}(f_1)>0$,  by replacing $g_k$ by $g^{-1}_k$ and switching the roles of $f_1$ and $f_2$, it is reduced to one of the sub-cases of case (II). Though it is possible to give an unified proof for both cases together, in the following we  give proofs for each of the cases. 

\medskip 
\noindent 
$\bullet$  Proof for case (I).

 In this case, we may assume that $v_m(h)>0$ for all the $L_k^p$-maps $h$ involved so that only the function $v_m$ is used. 

Now we choose $\epsilon_2$ for $U_{\epsilon_2}(f_2)$  as follows.    Since $f_2$ is $v_m$-stable, its  volume  $v_m(f_2)=\delta_2>0.$ Then there is a point $x_0\in {\bf CP}^n $ such that $v_m(df_2)(x_0)>0.$  Hence  there are   positive constants $\gamma$ and $\rho$ small enough such that for any
$x$ in the ball   $B(x_0; \rho)$  of radius $\rho$ centered at $x_0$,  $v_m(df_2)(x)>\gamma$. Then there is an ${\tilde \epsilon}_2>0$ such that for any $h$ with $\|h-f_2\|_{C^1}<{\tilde  \epsilon}_2$, the same is true. 

Now $\|h-f_2\|_{C^1}\leq C_2 \cdot \|h-f_2\|_{k, p}$ by our assumption.
Hence  we choose $\epsilon_2$ by the requirement that $\epsilon_2<{\tilde \epsilon}_2/C_2$ so that  for any $h\in U_{\epsilon_2}(f_2)$,  $\|h-f_2\|_{C^1}
<{\tilde \epsilon}_2$.

With this choice of $\epsilon_2$, for any $h\in U_{\epsilon_2}(f_2)$, and any point $x\in B(x_0; \rho)$,   
$v_m(dh)(x)>\gamma$.

In these notations, the condition (a) above implies that for $g_k$ in $G$
with $g_k=u_k\cdot D({a}_k)\cdot v_k$, we have     $lim_{k\mapsto \infty} a_k=0.$

After  taking subsequence, we may assume that $lim_{k\mapsto \infty} u_k=u$ and
$lim_{k\mapsto \infty} v_k=v$  in $SU(n+1)$

 Note that when considered as  automorphisms on ${\bf CP}^n$, the convergence  here are with respect to $C^{\infty}$-topology on the  corresponding mapping space.
 
 Let  ${\bf CP}^n={\bf C}^n\cup {\bf CP}^{n-1}  $. Now we  have the following  two cases:
 (A)  ${\bf CP}^{n-1} \cap  v(B(x_0; \rho) =\varphi $ and (B) ${\bf CP}^{n-1} \cap  v(B(x_0; \rho) \not =\varphi. $ Since $v$ is an  isometry,  it is easy to see that in both cases, there exits an $x_1\in B(x_0; \rho) $ and a positive number $\rho_1<< \rho$ such that $v(B(x_1; \rho_1))\cap {\bf CP}^{n-1}=\varphi.$
  We may assume that $dist (v(B(x_1; \rho_1)),{\bf CP}^{n-1})>\delta>0$. Then 
  for $i>i_0$ large enough, $dist (v_i(B(x_1; \rho_1)),{\bf CP}^{n-1})>\delta$ as well.
 Hence
 there is a large $R$ such that  for $i>i_0$ large enough, $v_i(B(x_1; \rho_1))$ is lying inside 
  $D^n_R=:D(R_1)\times \cdots \times D(R_n)\subset {\bf C}^n$  with $R_1=R_2\cdots =R_n=R$ of the $n$-fold  product of the  open disks centered at origin  in ${\bf C}$  with  radius $R$.

 Note that  in term of the   coordinate  of  ${\bf C}^n\subset {\bf CP}^n$, the action of $\Delta({a_i})$  is given by $\Delta({a_i})(z)= (a_i\cdot z_1, r_{i, 2} \cdot z_2, \cdots, r_{i, n} \cdot z_n)$ with $a_i\leq r_{i, 2}\cdots \leq r_{i, n}\leq 1.$

 Hence for any fixed $R>0 $ and any given $\epsilon >0, $    our assumption that $a_i\mapsto 0$ implies that there is a fixed  $i_0(\epsilon)>>0$ such that when $i>i_0(\epsilon),$
 $  \Delta({a_i})(v_i(B(x_1; \rho_1)))\subset a_i\cdot D(R_1)\times r_{i, 2}\cdot  D(R_2)\times \cdots \times r_{i, n}\cdot D(R_n)\subset D(\epsilon )\times D(R_2)\times \cdots \times D(R_n)$ with $R_k=R$ for $ k=1, \cdots, n.$

 Hence the 
 $vol ( \Delta({a_i})(v_i((B(x_1; \rho_1))))$ $\leq C_3\epsilon^2 R^{2n-2}$.  Here the  volumes  are computed with respect to the Fubini-Study metric which is uniformly equivalent
 to the flat metric on $D^n_R$ for fixed $R$. 
 Applying this to $g_i=u_i\circ \Delta({a_i})\circ v_i$, since $u_i$  preserves the  Fubini-Study metric, we conclude that
 for $i$ large enough,   $vol (g_i(B(x_1; \rho_1)  <C_4\epsilon^2.$
 
 Now
 
 $$v_m  ((h_i\circ g_i)|_{ B(x_1; \rho_1)})=v_m (h_i|_{ g_i( B(x_1; \rho_1))})$$  
$$ \leq ||h_i||^{2n}_{C^1} vol ( g_i( B(x_1; \rho_1)))\leq C_4\cdot ||h_i||^{2n}_{C^1}\epsilon^2$$
 $$  \leq C_4 ||h_i||^{2n}_{k, p} \epsilon^2  \leq  C_4(\|f_1\|_{k, p}+ ||h_i-f_1||_{k, p})^{2n}\epsilon^2 $$ $$ \leq C_4(\|f_1\|_{k, p}+\epsilon_1)^{2n} \epsilon^2 . $$
 By letting $\epsilon \mapsto 0,$ we conclude that  $\lim_{i\mapsto\infty}v_m ((h_i\circ g_i)|_{B(x_1; \rho_1)})
 =0.$
 
 Now  $h_i\circ g_i\in U_{\epsilon_2}(f_2)$. Recall that  for any $h\in U_{\epsilon_2}(f_2)$,   $v_m(h|_{ B(x_1; \rho_1)})>C\cdot \gamma\cdot \rho^{2n}_1$, which is   a fixed positive constant. 
 This is a contradiction.

\QED

\medskip 
\noindent 
$\bullet$  Proof for case (II).

The  proof for this case is a modification of the proof above for case (I). We still argue by contradiction. The choices for $\epsilon_1$ and $K_k$ are  the same as above.

Now we choose $\epsilon_2$ for $U_{\epsilon_2}(f_2)$  as follows.    Recall in this case, we may assume that $v_m(f_2)=0$ but $v_{m-1}(f_2)>0.$ Assume that  $v_{m-1}(f_2)=\delta_2>0.$ Then there is a point $x_0\in {\bf CP}^n $ such that $v_{m-1}(df_2|_U)(x_0)>0$ for some small open $K$-subset $U$.  Hence  there are   positive constants $\gamma$ and $\rho$ small enough such that (1) there is a $(m-1)$-dimensional geodesic ball   $B_{m-1}(x_0; \rho)$  of radius $\rho$ centered at $x_0$ and (2) for any
$x$ in  $B_{m-1}(x_0; \rho)$, the   $v_{m-1}(df_2|_{B_{m-1}(x_0; \rho)})(x)>\gamma$. Then there is an ${\tilde \epsilon}_2>0$ such that for any $h$ with $\|h-f_2\|_{C^1}<{\tilde  \epsilon}_2$, the same is true: $v_{m-1}(dh|_{B_{m-1}(x_0; \rho)})(x)>\gamma$. 

 As before, since  $\|h-f_2\|_{C^1}\leq C_2 \cdot \|h-f_2\|_{k, p}$ 
by  choosing  $\epsilon_2$  such  that $\epsilon_2<{\tilde \epsilon}_2/C_2$, we have  that  for any $h\in U_{\epsilon_2}(f_2)$,  $\|h-f_2\|_{C^1}
<{\tilde \epsilon}_2$.

With this choice of $\epsilon_2$, for any $h\in U_{\epsilon_2}(f_2)$, and any point $x\in B_{m-1}(x_0; \rho)$,   
$v_{m-1}(dh|_{B_{m-1}(x_0; \rho)})(x)>\gamma$. 

Let  ${\bf CP}^n={\bf C}^n\cup {\bf CP}^{n-1}  $.
Since $m-1=2n-1>\dim ({\bf CP}^{n-1})$, 
as before, we may assume that 
there exits an $x_1\in B_{m-1}(x_0; \rho) $ and a positive number $\rho_1<< \rho$ such that $v(B_{m-1}(x_1; \rho_1))\cap {\bf CP}^{n-1}=\varphi$ where $B_{m-1}(x_1; \rho_1)$ is the ($(m-1)$-dimensional) ball of radius $\rho_1$  centered at $x_1$  inside  $B_{m-1}(x_0; \rho)$.
 Then 
 $dist (v(B_{m-1}(x_1; \rho_1)),{\bf CP}^{n-1})>\delta>0$ and  
 $dist (v_i(B_{m-1}(x_1; \rho_1)),{\bf CP}^{n-1})>\delta$
for $i$ large enough. 

Hence
there is a large $R$ such that  for $i>i_0$ large enough, $v_i(B_{m-1}(x_1; \rho))$ is lying inside 
$D^n_R=:D(R_1)\times \cdots \times D(R_n)\subset {\bf C}^n\simeq {\bf R}^{2n}$  of the $n$-fold  product of the  open disks centered at origin  in ${\bf C}$  with  radius $R_i=R$.

Now the key point is to show that $vol_{m-1}[\Delta({a_i})(v_i(B_{m-1}(x_1; \rho_1)))]$ tends to zero as $i$ goes to infinity.
By projecting to one of the $(2n-1)$-dimensional coordinate planes, the tangent planes  of $v(B_{m-1}(x_1; \rho_1)))$ and $v_i(B_{m-1}(x_1; \rho_1)))$ at $x_1$ can be realized a graph of a linear function. Hence for $\rho_1$ small enough, $v(B_{m-1}(x_1; \rho_1)))$ and $v_i(B_{m-1}(x_1; \rho_1)))$ 
 can be realized a graph of a function as well over the coordinate plane with the dimension $2n-1=m-1$. Hence at least one of the first two coordinate lines in ${\bf R}^{2n}\simeq {\bf C}^n$ has to be included in the  $(2n-1)$-dimensional coordinate plane.

Now recall that
in term of the   coordinate  of  ${\bf C}^n\subset {\bf CP}^n$, the action of $\Delta({a_i})$  is given by $\Delta({a_i})(z)= (a_i\cdot z_1, r_{i, 2} \cdot z_2, \cdots, r_{i, n} \cdot z_n)$ with $a_i\leq r_{i, 2}\cdots \leq r_{i, n}\leq 1.$ By renaming  the coordinate line transversal to the  $(2n-1)$-dimensional coordinate plane as the last coordinate line of ${\bf R}^{2n}$, we may assume that  the  $(2n-1)$-dimensional coordinate plane are given by the coordinates $(x_1, \cdots, x_{2n-1})$. The  action $\Delta({a_i})$ then has the form
 $\Delta({a_i})(x_1,\cdots, x_{2n-1}, x_{2n})= (a_i\cdot x_1, r_{i, 2}  \cdot x_2, \cdots,  r_{i, K-1} \cdot x_{K-1},r_{i, K} \cdot x_{K} ,\cdots, r_{i, 2n} \cdot x_{2n})$ with the property that (1) $a_i\leq r_{i, 1}\cdots \leq r_{i, 2n-1}\leq 1;$ (2) there is a $K\geq 2$ such that 
$\lim_{i\mapsto \infty }r_{i, j}=0$ for $j<K$ and $r_{i, j}$ is bounded below from zero for $K\leq j\leq 2n-1;$ (3) $r_{i, 2n}\leq 1$.
Then $\Delta({a_i})$ is decomposed as $\Delta({a_i})=\Delta^{< K}({a_i})\circ \Delta^{\geq K}({a_i})$. Here the action of $\Delta^{< K}({a_i})$ is the same as that $\Delta({a_i})$ on the first $(K-1)$ variables but  the identity on  the rest of variables; and $\Delta^{\geq K}({a_i})$ is just another way around. 

Now using the family of $(2n-(K-1))$-dimensional planes parallel to the coordinate plane of the last $(2n-(K-1))$ coordinates of ${\bf R}^{2n}$ to slice $v(B_{m-1}(x_1; \rho_1)))$ and $v_i(B_{m-1}(x_1; \rho_1))),$ they are  realized as families of $(2n-K)$-dimensional balls (generically) of bounded volumes over the  $(K-1)$-dimensional finite ellipsoids $E$ and $E_i$   in the  coordinate plane of the first $(K-1)$ coordinates.  Denote these families  by 
$\pi: v(B_{m-1}(x_1; \rho_1)))\rightarrow E$ and $\pi_i: v_i(B_{m-1}(x_1; \rho_1)))\rightarrow E_i.$  Let $M$  be  an up bound of  $(2n-K)$-dimensional volume of the  balls (fibers )  in the family $\pi$. Since $v_i\mapsto v$ in $C^{\infty}$-topology, we may assume that $M$ is also  an up bound of the  $(2n-K)$-dimensional volumes  of the fibers  of   the family $\pi_i$ for $i$ large enough. Then by Fubini's theorem,  $vol_{m-1}(v_i(B_{m-1}(x_1; \rho_1)))\leq M\cdot  vol_{K-1}(E_i)$. Since $\Delta^{\geq K}({a_i})$ acts as identity in the first $K-1 $ variables with a dilation factor less than $1$ on the rest of variables, we still have 
$$vol_{m-1}(\Delta^{\geq K}({a_i})(v_i(B_{m-1}(x_1; \rho_1))))\leq M\cdot  vol_{K-1}(E_i).$$

Now applying $\Delta^{< K}({a_i})$ by noting that $\Delta^{< K}({a_i})$  only acts on the first $K-1$ variables leaving the rest unchanged, 
$$vol_{m-1}(\Delta({a_i})(v_i(B_{m-1}(x_1; \rho_1))))$$ $$ =vol_{m-1}(\Delta^{<K}({a_i})\circ \Delta^{\geq K}({a_i})(v_i(B_{m-1}(x_1; \rho_1))))$$ 
$$ \leq M\cdot  vol_{K-1}(\Delta^{< K}({a_i})(E_i)).$$ 

Now $E$ and hence $E_i$ are lying in a finite  ball of a large radius $R$ in  the coordinate plane of the first $K-1$ variables. 

\noindent
By the definition, $vol_{K-1}(\Delta^{< K}({a_i})(E_i))\leq r_{i, K-1}vol_{K-1}(E_i))\leq C\cdot r_{i, K-1}R^{K-1}$ that goes to zero as $i$ goes to infinity. Hence
$ vol_{m-1}(\Delta(a_i)(v_i(B_{m-1}(x_1; \rho_1))))$
 $\mapsto 0$ as $i\mapsto \infty$. Note that in  above the volume is computed in Euclidean metric. Since $v_i(B_{m-1}(x_1; \rho_1))$ and $\Delta(a_i)(v_i(B_{m-1}(x_1; \rho_1)))$ are lying in a  bounded  region $D_R^n\subset {\bf C}^n$, the same conclusion is true with respect to the Fubini-Study metric. Then as before, we conclude from above that $vol_{m-1} (g_i(B_{m-1}(x_1; \rho_1) \mapsto 0$ as $i\mapsto \infty$.

 The rest of the proof is the same as the one for case I by replacing $B(x_1; \rho_1)$  there by  $B_{m-1}(x_1; \rho_1)$ and  $m$-volumes
 there by the corresponding $(m-1)$-volumes. 
  
  Indeed
 $$v_{m-1} ((h_i\circ g_i)|_{ B_{m-1}(x_1; \rho_1)})=v_{m-1} (h_i|_{ g_i( B_{m-1}(x_1; \rho_1))})$$  
 $$ \leq ||h_i||^{2n-1}_{C^1} vol_{m-1} ( g_i( B_{m-1}(x_1; \rho_1)))  $$ $$\leq C(\|f_1\|_{k, p}+ ||h_i-f_1||_{k, p})^{2n-1} vol_{m-1} ( g_i( B_{m-1}(x_1; \rho_1)))  $$ $$ \leq C(\|f_1\|_{k, p}+\epsilon_1)^{2n-1}  vol_{m-1} ( g_i( B_{m-1}(x_1; \rho_1))). $$
 We conclude that  $\lim_{i\mapsto\infty}v_{m-1} ((h_i\circ g_i)|_{B_{m-1}(x_1; \rho_1)})
 =0.$
 
 Now  $h_i\circ g_i\in U_{\epsilon_2}(f_2)$. Recall that  for any $h\in U_{\epsilon_2}(f_2)$,   $v_{m-1}(h|_{ B(x_1; \rho_1)})>C\cdot \gamma\cdot \rho^{2n-1}_1$, which is   a fixed positive constant. 
 This is a contradiction.

\QED
\begin{cor}
	Corollary 1.1 holds.
\end{cor}
\proof 

 We need to show that  $\Gamma_f$ is compact  if $f$ is $v$-stable. 

By taking $f=f_1=f_2$ with $f$ being  $v$-stable,   the  above theorem implies that there is a compact subset $K\subset G$ such that $\Gamma_f$ is contained in $K$.
It is well-known that the action map $\Phi:G\times {{\cal M}}_{k,p}(M, N) \rightarrow {{\cal M}}_{k,p}(M, N)$ is continuous (for  a proof, see  \cite{Liu1} for instance )
 so that $\Gamma_f$ is closed.

\QED

\begin{cor}
	The Corollary  1.2 holds.
	
\end{cor}

\proof

We  need to  prove  that the  $G$-orbit  of $f$  in
$   { {\cal  M}}^{*}_{k, p}$ is closed.

Rename $f$ as  $f_1$. If the corollary is not true, there exist $g_i\in G$ and $f_2\in  { {\cal  M}}^{*}_{k, p} $ such that
$f_2=\lim_{i\mapsto \infty}f_1\circ g_i,$ but $f_2$ is not in $G\cdot f_1$.
Therefore for any $U_{\epsilon_2}(f_2), $  when $i$ is large enough,  $f_1\circ g_i$ is in $U_{\epsilon_2}(f_2). $
On the other hand,  the Theorem 3.1  with the same notation implies that for all such $i$,  $g_i$ is in the compact set $K_1$. Therefore, we may assume that
$\lim_{i\mapsto \infty} g_i=g$ in $K_1.$ Consequently, $f_2=\lim_{i\mapsto \infty}f_1\circ g_i=f_1\circ g.$ That is $f_2\in G\cdot f_1$ which is  a contradiction.

\QED

We restate the  Propostion 1.1 below.

\begin{pro}
	Let $G=PSL(n+1, {\bf C})$.	Then   the  space 
	 ${{\cal M}}^*_{k,p}({\bf P}^n, N)$ of  $v$-stable  $L_k^p$-maps  on ${\bf P}^n$  is $G$-Hausdorff in the sense that for any two diffent $G$-orbits $Gf_1$ and $Gf_2$, there exit $G$-neighborhoods $G U_1$ and $GU_2$ such that $GU_1 \cap G U_2=\varphi$. Therefore, the  quotient space  ${{\cal M}}_{k,p}({\bf P}^n, N)/G$  of unparametrized $v$-stable  $L_k^p$-maps is Hausdroff.
\end{pro}

\proof

By Theorem 3.1, for any $g$ not in the compact set $ K_1$ and $h\in  U_{\epsilon_1}(f_1)$, $h\circ g$ is not in $U_{\epsilon_2}(f_2)$. By our assumption, we may assume that
$  U_{\epsilon_1}(f_1)$ and $U_{\epsilon_2}(f_2)$ have no intersection.

\vspace{2mm}
\noindent
\vspace{2mm}
\noindent $\bullet$ Claim: when $\epsilon_i,i=1,2$  are small enough,
$(G\cdot U_{\epsilon_1}(f_1))\cap U_{\epsilon_2}(f_2)$ is empty.

\proof

\vspace{2mm}
\noindent
If this is not true, there  are $h_i\in U_{\delta_i}(f_1)$ and $g_i\in K_1$
such that $h_i\circ g_i$ is in $ U_{\delta_i}(f_2)$ with $\delta_i\mapsto 0.$
The compactness of $K_1$ implies that after taking a subsequence, we have that
$\lim_{i\mapsto \infty}g_i=g\in K_1.$ Since $\delta_i\mapsto 0$, we have that
$f_1=\lim_{i\mapsto \infty}h_i$ and $f_2=\lim_{i\mapsto \infty}h_i\circ g_i=f_1\circ g.$
Hence, $f_1$ and $f_2$ are in the same orbit which contradicts to our assumption.
Note that in the last identity above, we have used the fact that the action map
$\Psi: G\times  {{\cal M}}_{k,p}({\bf P}^n, N) \rightarrow {{\cal M}}_{k,p}({\bf P}^n, N)  $ is continuous.

\QED

Of course the same proof also implies that
$(G\cdot U_{\epsilon_2}(f_2))\cap U_{\epsilon_1}(f_1)$ is also empty for sufficiently
small $\epsilon_i, i=1, 2.$

If $h\in (G\cdot  U_{\epsilon_1}(f_1))\cap ( G\cdot  U_{\epsilon_2}(f_2)), $  then there are
$h_i\in U_{\epsilon_i}(f_i) $ and $g_i\in G, i=1, 2$ such that
$h=h_1\circ g_1=h_2\circ g_2. $ Hence $h_2=h_1\circ g_1\circ g_2^{-1}$ and
$(GU_{\epsilon_1}(f_1))\cap U_{\epsilon_2}(f_2)$ is not empty. This contradicts
to the above claim.

\QED

\subsection {{  Proof  of Theorem 1.2 and  Theorem 1.3}}

We need show that 	(i) any $v_{2n}$-stable  $L_k^p$-map  $f$ with the domain ${\bf CP}^n$  is stable; (ii) if $v_{2n}(f)=0$ but $v(f)>0$, then  $f$ is either a stable map  or  one of $S^1$-invariant maps (up to an conjugation).

\proof

Since in both cases $v_(f)>0$,  ${\Gamma}_{f}$ is compact.  If the identity component ${\Gamma}^0_{f}$  is  nontrivial,
let ${\tilde {\Gamma}}^0_{f}$ be the lifting of ${\Gamma}^0_{f}$ in $SL(n+1,{\bf C})$. Since 
maximal compact and  connected  subgroup  $SU(n+1)$ has only one orbit inside $SL(n+1,{\bf C})$ under conjugations, we may assume that ${\tilde {\Gamma}}^0_{f}$ is  contained in $SU(n+1)$ so that ${\Gamma}^0_{f}$ is contained in $SU(n+1)/{\bf Z}_2^{n+1}$. Thus up to a conjugation, 
 we may assume that  ${\tilde {\Gamma}}^0_{f}$ contains a subgroup   $S^1$ inside   the maximal tours   $T^{n+1}\subset SU(n+1)$  such that the induced  action of  $S^1\subset {\tilde {\Gamma}}^0_{f} $  on ${\bf P}^n$  given by $e^\theta\cdot (z_0:z_1:\cdots:z_n)=  (e^{m_0\theta}z_0:e^{m_1\theta}z_1:\cdots:e^{m_n\theta}z_n)$ is nontrivial.  It has the form $e^\theta\cdot (z_1,  \cdots, z_n)=  (e^{m_1\theta}z_1, e^{m_2\theta}z_2, \cdots, e^{m_n\theta}z_n)$ with some of $m_i\not = 0, $  with respect to  the coordinate $z=(z_1, \cdots,z_n)$ of ${\bf C}^n$ for a  proper choice of the decomposition     ${\bf CP}^n={\bf C}^n\cup {\bf CP}^{n-1}.$ 
  The $S^1$-action here is  the so called  (one of ) standard  $S^1$-action in (ii) above.

  If $\dim (L(\Gamma_f^0))\geq 2$, let $\xi_1$ and $\xi_2$ be two linear independent elements in $L(\Gamma_f^0)=L({\tilde {\Gamma}}^0_{f})$ with the corresponding $1$-parameter subgroups $S^1_{\xi_1}$ and $S^1_{\xi_2}$ both $\simeq S^1$,  each acting on ${\bf P}^n$ in the form above up to a conjugation. Since both action come from/extend to  the corresponding $C^*$-actions,  if the two actions of $S^1_{\xi_1}$ and $S^1_{\xi_2}$ are identical on some open set of ${\bf P}^n$, so are the corresponding  $C^*$-actions. Since 
   the  $C^*$-actions are algebraic,   the two actions of $S^1_{\xi_1}$ and $S^1_{\xi_2}$ are the same if they agree on an open set.
  
  Now consider the  map from $L(\Gamma_f^0)$ to the set of vector fields ${\cal V}({\bf P}^n)$ given by sending $\xi$ to the vector field
  $X_{\xi}$ that generates the action of $S^1_{\xi}.$
  The special form of the $S^1$-action above implies that the map  is injective. Let $L(\xi_1, \xi_2)$ be the linear span of $\xi_1 $ and $\xi_2$. Then above argument implies that the "plane" field $L(X_{\xi_1},X_{\xi_2})$ is $2$-dimensional generically in the sense that the set of points $x$ where $\dim[L(X_{\xi_1},X_{\xi_2})(x)]<2$ has no interior point. Now by definition $L(X_{\xi_1},X_{\xi_2})(x)\subset L(\Gamma_f^0)$ is  contained in $\ker (df_x)$. Hence at any  generic point $x$, the $rk (df_x)<m-1$. Hence  the same is true for all point $x\in {\bf P}^n$ by the lower semi-continuity of the function  $x\rightarrow rk (df_x)$, which contradicts to $v_{m-1}(f)>0.$ Hence $\dim  L(\Gamma_f^0)=1$ and $\Gamma_f^0\simeq S^1.$
  
   We have proved that if $v(f)>0$ and    
${\Gamma}^0_{f}$  is  nontrivial,  then $f$ is one of the  standard  $S^1$-invariant map up to a conjugation. This proves (ii).

To prove (i), still assume the identity component ${\Gamma}^0_{f}$  is  nontrivial, hence with the $S^1$-action above.
Assume that $m_1\not = 0. $ Since  $f$  is   $v_{2n}$-stable,   $df_v$ is injective at some point $x_0$ so that $f_v$ is  a local embedding on a small $\epsilon$-polydisk  
$D^n_{\epsilon }(x_0)=:D_{\epsilon}(x_{0, 1})\times \cdots \times D_{\epsilon}(x_{0, n}) \subset {\bf C}^n$ where $D_{\epsilon}(x_{0, j})$ is the   open disk in ${\bf C}$  centered  at $x_{0, j}$    with  radius $\epsilon$.

By taking a smaller $\epsilon$-polydisk, we may assume that the origin $0\in {\bf C}$, which is  the fixed point of   the  nontrivial   action $\phi_i(\theta, z)= e^{m_i\theta}z $  is not contained in $D_{\epsilon}(x_{0, i})$.

Then for $\phi\not={\it id}\in S^1$ but sufficiently close to ${\it id}$, the image  $\phi (D_{\epsilon/2}(x_{0, i}))$ is different from $D_{\epsilon/2}(x_{0, i})$ but still inside  $D_{\epsilon}(x_{0, i})$ if $\phi_i$ is nontrivial, where $\phi_i$ is the $i$-th component of $\phi$   so that the image  $\phi_i (D^n_{\epsilon/2}(x_{0}))$ is different from $D^n_{\epsilon/2}(x_{0}).$
  Since $f$ is an embedding on $D^n_{\epsilon}(x_0)$,  the images $f(D^n_{\epsilon/2}(x_0))$ and  $f(\phi(D^n_{\epsilon/2}(x_0)))=f\circ \phi(D^n_{\epsilon/2}(x_0))$ then  are different. Hence on $D^n_{\epsilon/2}(x_0)$, $f\circ \phi\not=f$ already so that $\phi$ is not lying inside   the isotropy group of $f$. This is a contradiction.

\QED

\end{large}

\end{document}